\documentclass[11pt]{article}%
\usepackage[margin=1in]{geometry}
\usepackage{mathpazo}
\usepackage{amsfonts}
\usepackage{amssymb}
\usepackage{graphicx}
\usepackage{setspace}
\usepackage[round]{natbib}
\usepackage{amsmath}%
\usepackage{amsthm}%
\usepackage{paralist}
\usepackage{nicefrac}
\usepackage{bm}
\usepackage{mathtools}

\usepackage{booktabs}
\usepackage{lineno}

\usepackage{bm}
\usepackage{float}
\usepackage{titlesec}
\usepackage{draftwatermark}

\SetWatermarkText{DRAFT}
\SetWatermarkScale{0.8}
\SetWatermarkColor[gray]{0.99}

\usepackage{epigraph}
\usepackage{nicefrac}

\usepackage[font={small,it}]{caption}

\usepackage{mathtools}
\usepackage{lscape}
\usepackage{subcaption}
\usepackage{suffix}
\usepackage{color}
\usepackage{setspace}
\usepackage{mathrsfs}
\usepackage[inline]{enumitem}
\usepackage{framed}
\usepackage{multirow}
\usepackage{wrapfig}

\usepackage{qtree}
\usepackage{algorithm,algorithmic}

\RequirePackage[hyphens]{url}
\usepackage{xcolor}
\usepackage{hyperref}
\definecolor{darkblue}{rgb}{0.0,0.0,0.3}
\hypersetup{colorlinks,breaklinks,linkcolor=darkblue,urlcolor=darkblue,
            anchorcolor=darkblue,citecolor=darkblue}
\usepackage{empheq}
\usepackage[most]{tcolorbox}

\setlist*[enumerate]{label=(\roman*)}

\def\boxit#1{\vbox{\hrule\hbox{\vrule\kern6pt
          \vbox{\kern6pt#1\kern6pt}\kern6pt\vrule}\hrule}}
					
\usepackage{tikz}
\usetikzlibrary{calc,shapes,bayesnet}

\newcommand{\e}{{\rm e}} % for exponentials

\newcommand{\y}{{\bf y}}

\newcommand{\E}{{\mathbb E}}

\newcommand{\I}{{\bf I}}

\newcommand{\X}{{\bf X}}

\newcommand{\bbeta}{\boldsymbol{\beta}}

\newcommand{\bmu}{\boldsymbol{\mu}}

\newcommand{\ben}{\begin{enumerate}}
\newcommand{\een}{\end{enumerate}}
\newcommand{\beq}{\begin{equation}}
\newcommand{\eeq}{\end{equation}}
\newcommand{\bde}{\begin{description}}
\newcommand{\ede}{\end{description}}

\newcommand{\vectornorm}[1]{\left|\left|#1\right|\right|}

\newcommand{\NormRV}{\mathcal{N}}

\newcommand{\dd}[2]{\frac{\partial #1}{\partial #2}}

\newcommand{\ind}{\stackrel{\mathrm{ind}}{\sim}}

\theoremstyle{slplain}

\newtheorem{theorem}{Theorem}

\newtheorem{remark}[theorem]{Remark}

\numberwithin{theorem}{section}

\singlespace

\title{Polynomial Log-Marginals and Tweedie's Formula \\ When Is Bayes Possible?\thanks{This preprint is a work in progress and may be updated; feedback and corrections are welcome.}}

\author{
  Jyotishka Datta \\ 
  Department of Statistics, Virginia Tech \\ 
  \texttt{jyotishka@vt.edu}
  \and
  Nicholas G. Polson \\
  Booth School of Business, University of Chicago \\
  \texttt{ngp@chicagobooth.edu}
}
\date{\today}

\begin{document}
\maketitle

\begin{abstract}
Motivated by Tweedie's formula for the Compound Decision problem \cite{robbins1951asymptotically}, we examine the theoretical foundations of empirical Bayes estimators that directly model the marginal density $m(y)$. Our main result shows that polynomial log-marginals of degree $k \ge 3 $ cannot arise from any valid prior distribution in exponential family models, while quadratic forms correspond exactly to Gaussian priors. This provides theoretical justification for why certain empirical Bayes decision rules, while practically useful, do not correspond to any formal Bayes procedures. We also strengthen the diagnostic by showing that a marginal is a Gaussian convolution only if it extends to a bounded solution of the heat equation in a neighborhood of the smoothing parameter, beyond the convexity of $c(y)=\tfrac12 y^2+\log m(y)$.
\end{abstract}


\section{Empirical Bayes and compound decision} 

Consider an experiment in which unknown parameters, say $\bmu = (\mu_1, \ldots, \mu_n)$, gives rise to observations $Y_i \ind f(y \mid \mu_i)$, $i = 1, 2, \ldots, n$, and the goal is to estimate $\bmu$ under a loss function $\ell(y, \mu)$, e.g., the squared-error loss $n^{-1}\vectornorm{\hat\bmu-\bmu}^2$. We assume that the $Y_i$'s are exchangeable, i.e., invariant under permutation of the indices, which implies conditional independence and identical nature by de Finetti's theorem. For Gaussian $f(\cdot)$, this problem is known as the normal means, or Gaussian sequence, or the Gaussian compound decision problem. In the compound decision problem \citep{robbins1951asymptotically}, the estimates (or any decision $\delta_i$ about $\mu_i$) is allowed to depend on all observations $\y = (y_1, \ldots, y_n)$, under the compound risk $\E_{\bmu}\{\ell(\delta, \bmu)\}$. The compound decision problem, together with Stein's shrinkage phenomenon \citep{stein1956inadmissibility}, shows that by allowing individual decisions/estimates to depend on the entire sequence, substantial reduction of total loss can be achieved. It is easy to see that the means $\mu_i$ are typically less extreme than the observations $y_i$; as Efron \citep{efron2011tweedie} emphasizes, such regression to the mean is essentially the selection-bias or ``winner's curse" phenomenon. The empirical Bayes (EB) methodology, introduced by \citep{robbins1956empirical}, provides a statistical procedure where the goal of approximating the ideal or oracle Bayes rule is nearly achieved without specifying a prior \citep{zhang2003compound, efron2011tweedie}. {A key appeal here is \emph{observability} using Tweedie's formula that depends only on the marginal $m(y)$, which is estimable from data without committing to a prior. Tweedie's formula is not only a convenient empirical Bayes identity but also sits at the core of Stein’s unbiased risk estimate (SURE), highlighting how shrinkage and risk estimation meet in the normal means model.} As noted by \citet{zhang2003compound}, the compound decision problem and the empirical Bayes methodology have been called the `two breakthroughs' in post-war statistics by \citep{neyman1962two}, as they have influenced a lot of modern statistical methods, especially in the context of high-dimensional data. {In practice, EB includes Type~II MLE/MMLE fits and related procedures that target sparse-signal regimes common in modern applications.}

One of the main attractions of the empirical Bayes framework for the compound decision problem is that it can result in risk reduction without solving the deconvolution problem of estimating the common conditional distribution $F$. Rather, one could work with the observable marginal density $m(y)$, and use either parametric or nonparametric models to estimate it, and use the estimated $\hat{m}(y)$ as a plug-in in the empirical Bayes approach. {This observable-only route computes estimators directly from $m(y)$, avoiding explicit priors and high-dimensional integration.} For example, assuming $\mu_i$'s to be normally distributed with zero mean and unknown variance $\sigma_0^2$, leads to a Gaussian marginal for the $Y_i$'s. Using a method of moments estimator for the $\sigma_0^2$ leads to the famous, classical James--Stein estimator $\delta(y) = (1 - (n-2)/\sum_{i=1}^{n} y_i^2) y$ \citep{stein1956inadmissibility}. Difficulty in choosing a parametric model for the marginal, leads the way to more flexible nonparametric models such as \citet{kiefer1956consistency}, which is easily seen to be convex and can benefit hugely from recent advances in convex optimization tools, as shown by \citet{koenker2014convex}. Yet another approach is proposed by Efron \citep{efron2008microarrays, efron2012large, efron2011tweedie}, using `Lindsey's method', which consists of estimating the marginal distribution from a Poisson regression of an integer order $K$, based on Lindsey's Method that models the histogram bin counts of the observed $\y$-values as Poisson random variables. \citet{efron2011tweedie}'s formulation yields locally adaptive shrinkage based on the local shape of the marginal histogram, while making no parametric assumptions on the prior. This has inspired subsequent developments in the empirical Bayes/compound decision problem by nonparametric density estimates of the marginal \citep[e.g.,][]{simon2013estimating, wager2014geometric}. As emphasized by \citet{efron2014two}, empirical Bayes admits two complementary strategies: \emph{f-modeling}, which models the marginal density $m(\cdot)$ on the observation ($y$) scale, and \emph{g-modeling}, which models the prior $F$ on the parameter ($\mu$) scale. The former have been relied upon heavily in applications, especially after \citet{robbins1956empirical}, while the latter have been predominant in the theoretical EB literature \citep{laird1978nonparametric, jiang2009general}.  {Beyond EB, the Kiefer--Wolfowitz NPMLE can be viewed directly as a compound-decision rule; ``no oracle needed'' results due to \citep{ritov2024no} show its Tweedie plug-in is minimax without reducing to an oracle Bayes benchmark.} {In high-dimensional prediction problems, it is also appealing to bypass the challenges of prior specification on a large parameter vector and instead work with the \emph{observable} one-dimensional marginal $m(y)$ via Tweedie's formula. Such a `predictive Bayes' viewpoint leverages $m(y)$ directly for calibrated shrinkage and prediction without specifying a high-dimensional prior \citep{padilla2018deconvolution}. In a recent discussion, \citet{efron2020pea} contrasts \emph{prediction} with \emph{attribution}, arguing that prediction is inherently easier and more forgiving than attribution, and is often `anti-parsimonious' and `a-probabilistic,' with theoretical efficiency replaced by empirical methods. He also shows that a two-step hierarchical approach—using Tweedie's formula for the local false discovery rate in wide data settings (large $p$, small $n$) can drastically reduce dimensionality while avoiding explicit prior modeling and remaining transparent about the observable target. This motivates our focus on which functional forms of $m(\cdot)$ can arise from a \textit{bona fide} Bayes model, and when the observable-only route, though appealing, might not always satisfy essential Bayes constraints.} 

Given the popularity and performance of EB methods, a natural question from a Bayesian perspective is: \textbf{for what class of prior distributions could a particular empirical-decision rule arise?} {Our focus is precisely on which functional forms of the \emph{observable} $m$ can actually occur as genuine Bayes marginals via Tweedie's identity.} In this article, we provide partial answers. First, if $\log m(y)$ is a polynomial of degree $>2$, no prior can produce $m(y)$ via Tweedie's formula; the only polynomial cases compatible with Bayes are degrees $\le 2$, with $K=2$ corresponding to a Gaussian prior. Second, beyond the convexity constraint on $c(y)=\tfrac12 y^2+\log m(y)$ as required by \citet{koenker2014convex}, we point out that a marginal is a Gaussian convolution if and only if it extends (in a neighborhood of the smoothing parameter) to a bounded solution of the heat equation $u_t=u_{xx}$ by the Hirschman-Widder Weierstrass representation theorem \citep[Thm.~VIII.6.3]{widder1955convolution}.

\subsection{Tweedie's formula}
In the simplest normal means model, let $\mu_i \stackrel{\text{iid}}{\sim} F$ and $Y_i \mid \mu_i \stackrel{\text{ind}}{\sim} \NormRV(\mu_i,1)$. A remarkable result is the Tweedie's formula\footnote{It seems the nomenclature for Tweedie's formula is an interesting case of Stigler's law of Eponymy. In his famous paper, \citet{robbins1956empirical} credits Tweedie, and \citet{efron2011tweedie} follows, but, in her blog, Jiaying Gu documents a chain from Laird--Louis to Tukey to \citet{dyson1926method}, and comments that `Tweedie's formula perhaps should be called Eddington's formula'. On the other hand, \citet{west1982aspects} provides the formulae for posterior moments for both location and scale and credits \citet{masreliez1975approximate}.}:
\begin{equation}
    \delta(y)=\E(\mu \mid Y=y)= y + s_F(y), \quad \text{where} \quad s_F(y)=\partial/\partial y \log m(y). \label{eq:tweedie-1}
\end{equation}
Here, $s_F(y)$ is the `Bayesian correction term' and $m(y)=\int \varphi(y-\mu)\,dF(\mu)$ is the marginal density of $y$. An attractive property of the Tweedie--Eddington formula \eqref{eq:tweedie-1} is that the Bayes decision $\delta(\cdot)$ is expressed as a function of the observations $\y$, not involving the unknown prior $F$ \citep{efron2011tweedie, ritov2024no}. 
\citet{efron2011tweedie} discusses the Tweedie's formula for the general exponential family. Start with a canonical parameter $\eta$ and cumulant generating function $\lambda(\cdot)$, and density at $y = 0$ being denoted by $m_0$, we have:
\[
y \mid \eta \sim f_{\eta}(y) = {\rm e}^{\eta y - \lambda(\eta)} m_0(y), \qquad \eta \sim g(\cdot).
\]
The posterior distribution under this model would be: $\pi(\eta \mid y) = f_{\eta}(y) g(\eta)/ m(y)$, where $m(y)$ is the marginal density given by: $m(y) = \int f_{\eta}(y) g(\eta) d \eta$, where the integral is carried over an appropriate sample space of the exponential family. \citet{robbins1956empirical} showed that the posterior density for $\eta$ given $y$ is also a member of the exponential family with canonical parameter given by the observation $y$ and cumulant generating function (CGF) $\lambda(y) = \log(m(y)/m_0(y))$. Differentiating the CGF yields the posterior mean and variance as:
\[
\E(\eta \mid y) = \lambda'(y), \quad {\rm Var}(\eta \mid y) = \lambda''(y).
\]
\citet{efron2011tweedie} also discusses how \citet{robbins1956empirical}'s famous empirical Bayes Poisson prediction formula $\E(\lambda \mid y) = (y+1) m(y+1)/m(y)$ can be recovered as an approximation of the posterior mean derived from the Tweedie's formula applied to the exponential family representation of the Poisson probability mass function. \citet{brown2013poisson} showed that this (now, more than) half-century old empirical Bayes formula can perform really well after some adjustments, viz., Rao-Blackwellization and smoothing (an isotonic regression to ensure montonicity). \citet{datta2016bayesian} propose a global-local shrinkage prior for the Poisson compound decision problem that can also account for inflation of small counts, a phenomenon called `quasi-sparsity'. 

\citet{polson1991representation} extends the Tweedie's formula to a general location model. Let $\hat{\mu}$ and $a$ be the MLE and the maximal ancillary respectively and $y = (\hat{\mu}, a)$. Then the posterior mean $\E(\mu \mid y)$ for the location model under an arbitrary likelihood $f$ and a normal prior $\mu \sim \NormRV(m, \tau^2)$ would be represented by the score function as:
\begin{equation}
    \E(\mu \mid y) = m - \tau^2 \dd{}{\hat{\mu}}\log p(\hat{\mu} \mid a), \label{eq:tweedie-polson}
\end{equation}
where $p(\hat \mu \mid a ) = \int p(\hat \mu \mid \mu, a) dF(\mu)$. \citet{polson1991representation} also shows that a similar score representation formula holds for a normal scale-mixture prior $\mu \sim \int \mathcal \varphi(\mu \mid m,\tau^2)\,dH(\tau^2)$, with $\tau^2$ in \eqref{eq:tweedie-polson} replaced by a suitable posterior quantity. \citet{polson2019bayesian} provides the Tweedie's formula for normal linear regression: 
\[
\y = \X\bbeta + \e, \qquad \e \sim \mathcal{N}(0,\Sigma).
\]
Let $g(\beta)$ denote the prior density of $\beta$, and define the marginal (prior predictive) density $m(y)=\int f(y \mid \beta)\,g(\beta)\,d\beta.$ Then, assuming $(X^\top\Sigma^{-1}X)^{-1}$ exists, the posterior mean of $\beta$ given $y$ is
\begin{equation}
    \mathbb{E}[\bbeta \mid \y]
= (\X^\top \Sigma^{-1} \X)^{-1}\, \X^\top \Big(\Sigma^{-1} \y + \nabla_{\y} \log m(\y)\Big). \label{eq:tweedie-reg}
\end{equation}
The gradient of the prior predictive score $ \nabla_{\y} \log m(\y)$ is the `Bayesian correction' term in the context of linear regression. Clearly, \eqref{eq:tweedie-reg} reduces to \eqref{eq:tweedie-1} for $\X = \I$, which reduces the regression problem to a normal means model. 

Finally, \citet{west1982aspects} provides us with the Tweedie's formula for the known location ($\mu = 0$), and unknown scale parameter $\sigma >0$, so that the likelihood is: $p(y \mid \sigma) = \sigma^{-1} p(\sigma^{-1}y)$. Define $\lambda = \sigma^{-2}$, then Theorem.~2.3.1 of \citet{west1982aspects} says: 

\begin{theorem}
Let the prior for the precision $\lambda>0$ be $\pi(\lambda)=\mathrm{Ga}\!\left(\tfrac{\alpha}{2},\tfrac{\beta}{2}\right)$ with density
$\pi(\lambda)\propto \lambda^{\alpha/2-1}\exp\!\{-\tfrac{\beta}{2}\lambda\}$, and let
\[
p_y(y)=\int_{0}^{\infty} \lambda^{1/2}\,p(\lambda^{1/2}y)\,\pi(\lambda)\,d\lambda,
\qquad 
g_y(y)=-\frac{\partial}{\partial y}\log p_y(y),\quad 
G_y(y)=\frac{\partial}{\partial y}g_y(y).
\]
Then the posterior moments of $\lambda$ are
\[
\mathbb{E}[\lambda\mid y]=\beta^{-1}\big\{(\alpha+1)-y\,g_y(y)\big\},
\qquad
\operatorname{Var}(\lambda\mid y)=\beta^{-2}\big\{2(\alpha+1)-3y\,g_y(y)-y^{2}G_y(y)\big\}.
\]
\end{theorem}

\vspace{0.2in}
\noindent \textbf{SURE:} Tweedie's formula is also related to the Stein's unbiased risk estimation \citep[SURE][]{stein1981estimation} framework for achieving a finite sample unbiased estimate of the prediction risk. For the normal means model, $\y \sim \NormRV(\bmu, \sigma^2 I)$, with predictions denoted by $\hat{\y}$, the SURE for total prediction error is: $\vectornorm{\hat{\y}-\y}^2 + 2 \sigma^2 \sum_{i=1}^{n}\dd{}{y_i}{\hat{y}_i}$. \citet{bhadra2019prediction} show that the SURE framework can be extended to regression by considering the singular value decomposition of the design matrix, and prove that the gain in predictive accuracy obtained by using the horseshoe prior can be achieved in finite samples, outperforming existing competitors, such as ridge regression and PCR. In a recent work, \citet{ghosh2025stein} show that taking $h(y)=s_F(y)$ in $\delta(y)= y + h(y)$ yields
\[
\mathrm{SURE}(\delta)=1+\frac{1}{n}\sum_{i=1}^n\Big\{ s_F(Y_i)^2 + 2\,\dd{}{y} s_F(Y_i) \Big\}.
\]
Choosing $F$ to minimize $\mathrm{SURE}(F)$ provides a data-driven approximation to the oracle Bayes rule. \citet{ghosh2025stein} further provide a unified `$g$-modeling' \citep{efron2014two} view showing that minimizing Stein's unbiased risk estimate is equivalent (up to a constant) to Hyv\'arinen's score matching \citep{hyvarinen2005estimation} for learning the prior, and establishes near-parametric convergence rates for regret/Fisher divergence with oracle inequalities under misspecification. For a comprehensive discussion of empirical Bayes and the Tweedie--Eddington formula, we refer readers to the lecture notes by Ignatiadis and Sen \citep{ignatiadis_sen_empirical_bayes_2025}.
 %\\ %

\section{Avoiding Deconvolution} 

All empirical Bayes methods, thus, share the goal of estimating the mixing measure $F$, irrespective of the loss function used \citep{gu2016problem}. The deconvolution problem, i.e., estimating $F$, was referred to as `estimating the inestimable' by Robbins \citep{koenker2024empirical}, and has a long history, as pointed out by \citet{ritov2024no}. However, the typical reason is to provide a solution to the normal means or compound decision problem \citep{robbins1956empirical}, i.e., a good estimate for the normal means $(\mu_1 , \ldots , \mu_n)$ something of interest in many different applications. For this problem, deconvolution is often an intermediate goal, but not the primary objective. 

We do not attempt a full review of the vast literature on the normal-means problem. Instead, we focus on a specific class of empirical-Bayes estimators called empirical decision rules \citep{koenker2014convex}. An empirical decision rule attempts to estimate the posterior mean function $\E(\mu_i \mid y_i) $ directly. Such a rule implicitly invokes the existence of some prior distribution $f_0$ but does not attempt to estimate it \citep{efron2011tweedie}. 

To understand this empirical-Bayes approach, return to the simple normal-mean model $(y_i \mid \mu_i) \sim \NormRV(\mu_i,1)$. Absent any prior information that would distinguish the $\mu_i $, a natural assumption is that they arise exchangeably from some prior, $ \mu_i \sim F$. An important result discussed by \citet{robbins1956empirical} holds that the posterior mean can be written as the unbiased estimate plus a Bayes correction:
\begin{equation}
    E(\mu|y)=y + \dd{}{y} \log m(y), \label{eq:tweedie}
\end{equation}
where $m(y) = \int \varphi(y \mid \mu) dF(\mu)$ is the marginal density of the data under this prior. The prior does not appear explicitly, but its effect is incorporated in $m(y)$. \citet{efron2011tweedie} nicely summarizes its appeal: ``{\textit{The crucial advantage of Tweedie's formula is that it works directly with the marginal density},}" thereby offering the statistician the opportunity to avoid deconvolution entirely. 

% \vspace{0.2in} \noindent \textbf
\subsection{Nonparametric Maximum Likelihood} 

Among methods that directly estimate the marginal, \citet{brown2009nonparametric} advocate for replacing $m(y)$ by a kernel density estimate of it, while the usage of a nonparametric maximum likelihood estimator \citep[e.g.][NPMLE]{kiefer1956consistency} for $m(y)$ has been investigated by various authors including \citet{jiang2009general, saha2020nonparametric, greenshtein2022generalized}. \citet{koenker2014convex} argue that recent advances in convex optimization has paved the way for wider usage and applicability of the Kiefer--Wolfowitz algorithm for the nonparametric empirical Bayes compound decision problems. Ritov \citep{ritov2024no} shows that the Tweedie plug-in based on the nonparametric MLE of the marginal (equivalently, the Kiefer-Wolfowitz NPMLE for the mixing law) yields a decision rule (based on the Tweedie's formula) that is minimax in both the empirical Bayes and compound-decision formulations, establishes concentration for the NPMLE decision, and links its risk directly to SURE - thereby providing a justification for NPMLE without oracle arguments. 

In other words, the key insight here is that, since we observe $m(y)$ directly, a non-parametric approach to estimate $m(y)$ and use that as a plug-in in Tweedie's formula is a natural solution. Moreover, since \eqref{eq:tweedie} only involves the marginal $m(y)$, one can bypass priors entirely and estimate: 
\[
\hat\mu_i^{\text{TW}} \;=\; y_i + \partial_y \log \hat m(y_i),
\]
where $\hat m$ is any nonparametric estimate of $m$ (e.g., a kernel density, a Kiefer--Wolfowitz NPMLE Gaussian mixture, or a shape-constrained fit), thus avoiding deconvolution and explicit prior modeling \citep{brown2009nonparametric, kiefer1956consistency, koenker2014convex, jiang2009general, saha2020nonparametric, greenshtein2022generalized}. Writing $\ell(y)=\log m(y)$ and $s(y)=\ell'(y)$, one may estimate $\ell$ by minimizing a score-matching/SURE objective with a smoothness penalty,
\[
\min_{\ell}\;\Big\{\;\frac{1}{n}\sum_{i=1}^{n}\big[s(y_i)^2+2\,s'(y_i)\big] \;+\; \rho\,\mathcal R(\ell)\;\Big\},
\qquad \mathcal R(\ell)=\int\big(\ell''(y)\big)^2\,dy,
\]
optionally under shape constraints such as convexity of $c(y)=\tfrac12 y^2+\ell(y)$ on a grid (linear inequalities), which yields a convex quadratic program and directly returns the Tweedie plug-in $\hat\mu_i^{\text{TW}}$ \citep[see also the SURE/score-matching equivalence in][]{ghosh2025stein}.

\subsection{Lindsey's method} \citet{efron2011tweedie} argues that the Tweedie's formula $y + \dd{}{y} \log m(y)$ requires a smoothly differentiable estimate of the log marginal $\log m(y)$, and describes a Poisson regression technique, called the `Lindsey's method' \citep{efron2004estimation, efron2011tweedie} and \citep[][ch.5]{efron2012large}. The Lindsey's method assumes that the log-marginal is a $K$-th degree polynomial, for some integer $K$. That is, it estimates the $\mu_i $  by first estimating the marginal density directly using a polynomial or spline on the log scale, e.g.,
\begin{equation}
    \hat{m} ( y) = \exp \left ( \hat{\beta}_0 + \sum_{k=1}^K \hat{\beta}_k y^k  \right ). 
\end{equation}
The $\beta_k$'s are estimated via a Poisson regression (using usual generalized linear model tools) to the histogram counts after binning the data. Under this form of the marginal density, the posterior mean function is easily computed by taking the log and differentiating. Following on the work of \citet{efron2011tweedie}, \citet{koenker2014convex} considered an estimator based on the observation that the posterior mean function $\E(\mu \mid y) $ in Tweedie's formula \eqref{eq:tweedie} is non-decreasing in $y$. This implies that the function
\[
c(y)= \frac{1}{2} y^2 + \log m(y) 
\]
is convex, a result which holds not merely for the Gaussian case, but for any exponential-family observational model, regardless of the prior $F$. However, an unconstrained or traditional kernel-based estimate of $m(y)$ does not ensure convexity of $c(y)$.  As discussed before, \citet{koenker2014convex} proposed to estimate $m(y)$ non-parametrically.
 
\section{A Bayesian view of empirical decision rules}  

From a Bayesian perspective, a natural question is: \textbf{for what class of prior distributions could a particular empirical-decision rule arise?} In the case of Efron's estimator, we provide an explicit answer in the following theorem: \textbf{there are no such priors, save the trivial case of a Gaussian prior} (and a quadratic log marginal). We state this as the following theorem: 

% \begin{theorem}\label{thm:polylog}
% Let $ \phi(y| \mu) = e^{ \mu y - h( \mu) } f_0(y ) $ be a linear exponential family with scalar parameter $\mu$ through the mixing measure $f_0$. If $K \geq 3$, then for all prior measures $F$ supported on $\Re$, and for all nonzero choices of  $ \beta_1 , \ldots , \beta_K $ then
% $$
% \int_{- \infty}^\infty \phi( y| \mu ) d G(\mu )  \neq \exp \left ( \hat{\beta}_0 + \sum_{k=1}^K \hat{\beta}_k y^k  \right )
% $$
% Hence, Efron's empirical-Bayes rule is therefore not a Bayes rule, in that it cannot arise from a valid application of Tweedie's formula under any prior. The only exception is $K = 2$, implying that $F$ is Gaussian.
% \end{theorem}

\begin{theorem}\label{thm:polylog}
Let $\varphi(y\mid\mu)=\exp\{\mu y-\tfrac12\mu^2\}\,\varphi(y)$ be the $\NormRV(\mu,1)$ exponential-family form, where $\varphi$ is the $\NormRV(0,1)$ density, and let $m(y)=\int \varphi(y\mid\mu)\,dF(\mu)$ be the marginal induced by a prior $F$ on $\mu\in\mathbb{R}$. If $\log m(y)$ is a polynomial of degree $K\ge 3$ on $\mathbb{R}$, then there is no such prior $F$. In other words, for all prior measures $F$ supported on $\mathbb{R}$, and for all nonzero choices of  $\beta_1 , \ldots , \beta_K $,
$$
\int_{- \infty}^\infty \phi( y \mid \mu ) ~d F(\mu )  \neq \exp \left ( {\beta}_0 + \sum_{k=1}^K {\beta}_k y^k  \right ), \quad \text{if} \quad K \geq 3.
$$
If $K=2$, then $F$ must be (possibly degenerate) Gaussian.
\end{theorem}

Hence, Efron's empirical Bayes rule is therefore not a formal Bayes rule, in that it cannot arise from a valid application of Tweedie's formula under any prior. The only exception is $K = 2$, implying that $F$ is Gaussian. Interestingly, when $\log m(y)$ is quadratic ($K=2$), Tweedie's formula gives $\E(\mu_i\mid Y_i = y)=(1-1/V)y$ for a normal marginal $\NormRV(0,V)$, and replacing $1/V$ by the unbiased estimator $(n-2)/\sum_{j=1}^n y_j^2$ yields the James--Stein estimator \citep{efron2011tweedie}. We note that the argument extends trivially to other smooth bases, such as spline expansions, since these do not define entire functions on the complex plane and therefore cannot yield a valid prior either. 

\citet{koenker2014convex} cast two EB estimators as convex programs: (i) a direct MLE of the marginal density $m$ under the shape constraint that $c(y)=\tfrac12 y^2+\log m(y)$ is convex : this enforces monotonicity of $\delta(y)=y+m'(y)/m(y)=c'(y)$ and yields a piecewise-linear $c$, so $\delta(y)=c'(y)$ is piecewise constant; and (ii) a Kiefer--Wolfowitz NPMLE of the mixing law $F$ via a finite-dimensional convex dual whose solution is atomic (with at most $n$ support points), with interior-point implementations delivering substantial speedups over EM while retaining excellent risk performance. For the \citet{koenker2014convex} estimator, all Gaussian convolutions do satisfy the convexity of $c$, but convexity is only necessary, not sufficient. Additional conditions are given by the following result adapted from Theorem~VIII.6.3 of \citet{widder1955convolution}.

We note that \citet{koenker2014convex} themselves emphasize the practical advantages of the weaker convexity constraint, viewing it as a deliberate relaxation rather than a deficiency\footnote{``\textit{(S)ince the monotonized Bayes rule imposes a weaker form of convexity constraint, it may have some advantages in misspecified settings where the original Gaussian location mixture assumption fails to be satisfied.}" \citep{koenker2014convex}}. Theorem \ref{th:widder} identifies precisely when such convexity-based rules cease to correspond to any Bayesian model.

\begin{theorem}\label{th:widder}[\citet{widder1955convolution}]
Let $ u: \mathbb{R} \times(0,\infty) \mapsto \mathbb{R}$ be $C^{2}$ in $x$, $C^{1}$ in $t$, and bounded for each $t>0$. 
Then the following are equivalent:
\begin{enumerate}
\item $u$ solves the heat equation $u_t=u_{xx}$ on $\mathbb{R}\times(0,\infty)$ and, for some bounded Borel measure $\mu$, 
\[
u(x,t)=\frac{1}{\sqrt{4\pi t}}\int_{\mathbb{R}} e^{-(x-y)^2/(4t)}\,d\mu(y)\qquad(t>0).
\]
\item For each fixed $t>0$, $u(\cdot,t)$ is the Weierstrass transform (Gaussian convolution) of a bounded function; in particular $m_t(x):=u(x,t)$ is a Gaussian convolution.
\end{enumerate}
If, in addition, $\mu$ has a bounded density $f$, then $u(\cdot,t)=\phi_t*f$ with $\phi_t$ the $N(0,2t)$ kernel. Consequently, for a single marginal $m(x)=u(x,t_0)$ (e.g., $t_0=1$), being a Gaussian convolution is stronger than having $c(x):=\tfrac12 x^2+\log m(x)$ convex: $m(\cdot)$ must extend to a bounded solution of the heat equation in a time neighborhood of $t_0$.
\end{theorem}

% \noindent\fbox{\parbox{0.98\linewidth}{
% \textbf{When can a Tweedie plug-in be Bayes?} 
% (i) If $\log m$ is a polynomial of degree $>2$, \emph{never} (Thm.~2). 
% (ii) If $\log m$ is quadratic, then $F$ is Gaussian (Thm.~2). 
% (iii) Beyond convexity of $c(y)=\tfrac12 y^2+\log m(y)$, $m$ is a Gaussian convolution only if it locally extends to a bounded solution of $u_t=u_{xx}$ (Thm.~3).
% }}

\begin{remark}{\textbf{One-parameter antecedents.}} For one-parameter models, early work showed important rigidity phenomena: for a class where the sample total is sufficient for any sample size, the fiducial distribution of the parameter cannot coincide with any Bayesian posterior under any prior \citep{grundy1956fiducial}. In one-parameter natural exponential families, \citet{sampson1975characterizing} shows that the moment generating function (as a function of the natural parameter) uniquely characterizes the family and provides necessary and sufficient conditions for such MGFs. In one-parameter settings where the posterior expectation of the population mean is a linear function of the sample observation, \citet{goldstein1975uniqueness} proves that the moments of the prior distribution are uniquely determined, giving precise uniqueness relations for such linear Bayes rules. Finally, in a one-parameter exponential family under squared-error loss, \citet{goldstein1977contractions} shows that the only Bayes estimators whose 1\text{-}Lipschitz (contraction) transforms are also Bayes are the linear (affine) estimators in the canonical sufficient statistic.
\end{remark}

%%% JD's addition

% \section{JD's addition}

% \vspace{0.5in}
% \noindent {\bf I have added the next part but the math needs to be checked for correctness etc. - JD}

\subsection{Proof of Theorem \ref{thm:polylog}}

\begin{proof}

% The proof is outlined below step-by-step:

% \begin{description}
% \item[Step 0 (Positivity).] 
First, note that for any probability prior $F$, the mixture $m(y)=\int\varphi(y-\mu)\,dF(\mu)$ is strictly positive on $\mathbb{R}$; thus $\log m(\cdot)$ is well-defined and finite on $\mathbb{R}$.

% \item[Step 1 (Gaussian tilt and entire mgf).]
Next, writing
\[
m(y)=\varphi(y)\int e^{\mu y-\mu^2/2}\,dF(\mu)=\varphi(y)\,M_H(y),
\]
where $H$ is the finite positive measure $dH(\mu)=e^{-\mu^2/2}\,dF(\mu)$ and $M_H(z)=\int e^{z\mu}\,dH(\mu)$ is its moment generating function. Then, for $z=x+it\in\mathbb{C}$,
\[
|M_H(z)|
\le \int e^{x\mu}e^{-\mu^2/2}\,dF(\mu)
= \int \exp\!\Bigl(-\tfrac12(\mu-x)^2+\tfrac12 x^2\Bigr)dF(\mu)
\le e^{x^2/2}.
\]
Hence $M_H$ is entire and of at most Gaussian growth in $\Re z$. Put $M:=H(\mathbb{R})=\int e^{-\mu^2/2}dF(\mu)\in(0,1]$ \footnote{Since \(F\) is a probability measure and \(0<e^{-\mu^{2}/2}\le 1\) for all \(\mu\in\mathbb{R}\), we have: $M:=H(\mathbb{R})=\int_{\mathbb{R}} e^{-\mu^{2}/2}\,dF(\mu)\in(0,1].$}

% On each strip \(|\Re z|\le R\),
% \[
% \bigl|e^{z\mu}e^{-\mu^{2}/2}\bigr|
% \le e^{R|\mu|-\mu^{2}/2}
% =\exp\!\Big(-\tfrac12\big(|\mu|-R\big)^{2}+\tfrac12R^{2}\Big)
% \le e^{R^{2}/2},
% \]
% so the integrand is dominated by an \(L^{1}(H)\) function independent of \(z\); hence \(M_H\) is entire by dominated convergence (and differentiation under the integral is justified).

% \item[Step 2 (Identity theorem $\Rightarrow$ Exponential of a polynomial)]
If $\log m(y)$ is a polynomial on $\mathbb{R}$, then
\[
\log M_H(y)=\log m(y)-\log\varphi(y)=:P(y)
\]
is also a polynomial (since, $-\log\varphi(y)=\tfrac12 y^2+\mathrm{const}$ is quadratic).
Both $M_H(z)$ and $e^{P(z)}$ are entire and agree for all real $z$; by the identity theorem for entire functions, $M_H(z)\equiv e^{P(z)}$ on $\mathbb{C}$ (see, e.g., \citealp{Ahlfors1979}). Evaluating on the imaginary axis and normalizing yields a characteristic function:
\[
\phi_{\widetilde H}(t):=\frac{M_H(it)}{M}=\exp\!\bigl(P(it)-\log M\bigr)=:\exp Q(t),
\]
so $\log\phi_{\widetilde H}$ is a polynomial.

% \item[Step 3 (Marcinkiewicz $\Rightarrow$ degree $\le 2$).]
By Marcinkiewicz's theorem, if a characteristic function has the form $\phi(t)=\exp\{Q(t)\}$ with $Q$ a polynomial on $\mathbb{R}$, then $\deg Q\le 2$. Moreover, $\deg Q=0,1$ correspond to a degenerate distribution and $\deg Q=2$ to a (nondegenerate) Gaussian distribution \citep[Thm.~2.5.3]{Bryc1995}; see also \citet[Ch.~3]{Lukacs1970} and \citet{Marcinkiewicz1939}.
Therefore $\deg P\le 2$, so $\deg\log m(\cdot)\le 2$. In particular, no prior $F$ can produce a marginal with polynomial $\log m$ of degree $K\ge 3$. \emph{(Marcinkiewicz's theorem applies here since $\phi_{\tilde H}(t)=M_H(it)/M$ is a bona fide characteristic function.)}

% \item[Step 4 (Quadratic case $\Rightarrow$ Gaussian $F$).]
If $\log M_H(y)=\alpha y^2+\beta y+\gamma$, then
\[
m(y)=\varphi(y)\,e^{\alpha y^2+\beta y+\gamma}
= C\exp\!\Bigl(-\tfrac12(1-2\alpha)\,y^2+\beta y\Bigr),
\]
with $C=e^\gamma/\sqrt{2\pi}$. Since $m$ is integrable, $1-2\alpha>0$, so $m$ is a (nondegenerate) Gaussian density.
Write $Y=\mu+\varepsilon$ with $\varepsilon\sim N(0,1)$ independent of $\mu\sim F$. Then $Y$ is normal and $\varepsilon$ is nondegenerate normal. By Cram\'er's decomposition theorem, $\mu$ must be (possibly degenerate) normal \citep[Thm.~2.5.2]{Bryc1995}; see also \citet{Cramer1936}. Hence, $F$ is Gaussian.
% \end{description}

\end{proof}

\begin{remark}
The argument extends verbatim to any LEF (linear exponential family) with quadratic carrier log-density, since Step~2 only uses that $\log f_0(y)$ is quadratic to pass from $\log m$ polynomial to $\log M_H$ polynomial. For general LEFs this implication is even more restrictive.
Background on characteristic functions and cumulants (used implicitly in Marcinkiewicz/Cramer arguments) is collected in \citet[Chap.~XV]{Feller1971}.
\end{remark}

% \section{Conclusion}
% Why do we care about whether an empirical Bayes estimator corresponds to a formal Bayes procedure, i.e., whether it can arise from a hierarchical model with an actual prior? One key reason is that the formal Bayes rules provide structural guarantees, such as admissibilty (or coherence or least-favorable minimaxity). For example, a standard result in Lehmann--Casella is the following. let $R(\theta,\delta)=\mathbb{E}_\theta L(\theta,\delta(X))$ be the risk and $r(\delta,\Pi)=\int R(\theta,\delta)\,d\Pi(\theta)$ the Bayes risk for a \emph{proper} prior $\Pi$. If $\delta_\Pi$ is Bayes with respect to $\Pi$ and $r(\delta_\Pi,\Pi)<\infty$, then $\delta_\Pi$ is admissible. Admissibility, itself, is usually not the end goal and one can easily construct `admissible' estimators that are not at all usable. As \citet{strawderman2021charles} puts it: ``Admissibility is a relatively weak optimality property, and the fact that a procedure is admissible does not give strong evidence that it should be adopted. Typically, there are a very large number of admissible procedures in any given problem, including all procedures which are unique Bayes with respect to any proper prior." Yet, guarantees in the form of such complete class theorems for `Bayes' decision rules are the price to pay to avoid doing full Bayes. We hope to have shown that while the EB/NPMLE plug-ins optimize observable criteria (SURE, regret), they will not automatically inherit these guarantees.

\section{Conclusion}
Why care whether an empirical Bayes estimator corresponds to a formal Bayes procedure, that is, whether it can arise from a hierarchical model with an actual prior? Proper Bayes rules bring structural guarantees: admissibility under a proper prior with finite Bayes risk, coherence of the posterior, and, in favorable cases, least-favorable-prior minimaxity. A standard result (see Lehmann and Casella) states: let $R(\theta,\delta)=\E_\theta L(\theta,\delta(X))$ be the risk and $r(\delta,\Pi)=\int R(\theta,\delta)\,d\Pi(\theta)$ the Bayes risk for a proper prior $\Pi$; if $\delta_\Pi$ is Bayes with respect to $\Pi$ and $r(\delta_\Pi,\Pi)<\infty$, then $\delta_\Pi$ is admissible. Admissibility itself is a weak filter and does not, by itself, recommend a procedure: as \citet{strawderman2021charles} notes, ``\textit{Admissibility is a relatively weak optimality property, and the fact that a procedure is admissible does not give strong evidence that it should be adopted. Typically, there are a very large number of admissible procedures in any given problem, including all procedures which are unique Bayes with respect to any proper prior.}" Our results mark the boundary for Tweedie-based rules: polynomial $\log m$ of degree greater than two cannot be Bayes, the quadratic case corresponds to a Gaussian prior, and a marginal is a Gaussian convolution only if it extends locally in the smoothing parameter to a bounded solution of the heat equation. We hope to have shown that empirical Bayes and NPMLE plug-ins optimize observable criteria such as SURE and regret, but they do not automatically inherit the guarantees that come with a proper Bayes formulation. %In short, proper Bayes delivers admissibility and, in favorable cases, least-favorable minimaxity, whereas Tweedie-NPMLE plug-ins trade those structural guarantees for observable-risk control (SURE/regret) and convex implementability.

\section*{Acknowledgments}
The first author (JD) gratefully acknowledges support from the National Science Foundation (NSF CAREER Award DMS-2443282).

% \section*{Acknowledgments}
% The first author (JD) gratefully acknowledges support from the National Science Foundation (NSF CAREER Award DMS-2443282).

\bibliographystyle{biometrika}
\bibliography{paper-ref}

\end{document}